\documentclass[reqno]{amsart}
\usepackage{amssymb,amsthm,amsfonts,amsmath,amscd,amssymb,epsfig}
\newcommand*{\rom}[1]{\expandafter\@slowromancap\romannumeral #1@}
\usepackage{enumerate,array,graphicx,color,amsthm,bbm,mathrsfs, a4wide}
\usepackage{spectralsequences}
\usepackage{rotating}
\usepackage{bbm}
\usepackage[all,cmtip]{xy}

\usepackage[latin1]{inputenc}
\usepackage{dsfont}
\usepackage[onehalfspacing]{setspace}%

\usepackage[CJKbookmarks,bookmarks=true,colorlinks=true]{hyperref}%
\usepackage[english]{babel}

\usepackage{mathabx}
\usepackage{mathtools}
\usepackage{indentfirst}
\usepackage{graphicx}
\usepackage{breqn}
\usepackage{verbatim}

\theoremstyle{plain}
\newtheorem{theorem}{Theorem}[section]

\newtheorem{lemma}[theorem]{Lemma}

\newtheorem{prop}[theorem]{Proposition}
\theoremstyle{definition}
\newtheorem{definition}[theorem]{Definition}

\newtheorem{example}[theorem]{Example}

\theoremstyle{remark}
\newtheorem{remark}[theorem]{Remark}

\numberwithin{equation}{section}

\makeatletter
\newcommand{\colim@}[2]{%
  \vtop{\m@th\ialign{##\cr
    \hfil$#1\operator@font colim$\hfil\cr
    \noalign{\nointerlineskip\kern1.5\ex@}#2\cr
    \noalign{\nointerlineskip\kern-\ex@}\cr}}%
}
\newcommand{\colim}{%
  \mathop{\mathpalette\colim@{\rightarrowfill@\textstyle}}\nmlimits@
}
\makeatother

\newcommand{\C}{{\mathbb{C}}}
\newcommand{\R}{{\mathbb{R}}}

\newcommand{\Z}{{\mathbb{Z}}}

\renewcommand{\epsilon}{\varepsilon}
\renewcommand{\phi}{\varphi}

\newcommand{\rvline}{\hspace*{-\arraycolsep}\vline\hspace*{-\arraycolsep}}
\newcommand{\bigzero}{\mbox{\normalfont\bfseries 0}}

\begin{document}
\title[Symplectic Homology and 3-dimensional Besse Manifolds with $c_1=0$]{Symplectic Homology and 3-dimensional Besse Manifolds with Vanishing First Chern Class}
\author{Do-Hyung Kim}
\date{\today}

\begin{abstract}
In this paper, we will show that certain types of symplectic homology can be used as an invariant of 3-dimensional Besse manifolds, which are strict contact manifolds with periodic Reeb flow. For simplicity, we will assume our Besse structures to be a trivial plane bundle. To identify Besse manifolds with such a condition, we actually compute the first Chern class of each Besse structure and classify the Besse manifolds with vanishing first Chern class. We will also compute Robbin-Salamon indices of periodic Reeb orbits in Besse manifolds, and symplectic homology (of its filling) of certain cases. From its definition, Besse manifolds naturally become Seifert fibration, and thus one can extract invariants such as the orbifold Euler characteristic and the Euler number of this Seifert fibration. These invariants will become important in our computation. 
\end{abstract}

\maketitle

\section{Background}

From its first appearance in \cite{SH0} and its generalization in \cite{SH1}, \emph{symplectic homology} has been a useful tool in symplectic/contact geometry and topology. It is a Floer-type homology defined on a special symplectic manifold with a contact boundary called a Liouville domain. Here, ``Floer-type" means that its differential is defined by cylinders determined by the Floer equation. While defining this object, one can choose admissible Hamiltonians that depend only on the radial coordinate. This makes the periodic Hamiltonian orbits, generators of the Floer chain complex, become periodic Reeb orbits on each slice with a fixed radial coordinate. By incorporating objects and properties in contact dynamics, symplectic homology relates the symplectic geometry of the Liouville domain and the contact geometry of its boundary.

Because of this nature, symplectic homology can be used as a tool for investigating properties in both contact and symplectic geometry. It is invariant under exact symplectomorphisms and symplectic deformations \cite{biased,Benedetti_2020}, so one can use symplectic homology to distinguish Liouville structures up to these equivalence relations. For example, \cite{mcleanstein} used symplectic homology to distinguish Stein structures, a special type of Liouville structure, on $\mathbb{R}^{2n}$ up to exact symplectomorphisms and symplectic deformations along Stein structures. Symplectic homology can also be used to prove the non-fillability property of particular contact manifolds, using the simple fact that symplectic homology can always be defined on an exact filling of a contact manifold. Zhengyi Zhou \cite{Zhou_2021,Zhou_2024} utilized this fact and showed that contradictions occur due to topological reasons, hence deducing the non-existence of exact fillings of certain contact manifolds. One can utilize other versions of symplectic homology, such as positive symplectic homology, (positive) $S^1$-equivariant symplectic homology, Rabinowitz-Floer homology, etc., or modify it by truncating the homology according to the action value to explore the geometry or dynamics in a more detailed fashion.

In this paper, we employ some variants of symplectic homology as invariants among a specific type of contact manifolds. This means that for any contact manifold of this type, we first choose a Liouville filling of the contact manifold and then compute the symplectic homology of the filling. Since symplectic homology generally depends on the filling rather than the contact boundary itself, we must also ensure that the symplectic homology is independent of our choice of filling in order to properly use it as an invariant.

The class of contact manifolds we are going to work with is called \emph{Besse manifolds}, which are strict contact manifolds (i.e., contact manifolds with a given contact form) such that all of their Reeb orbits are periodic. Besse manifolds have rich symplectic/contact geometric and dynamical properties. For instance, there is a structure theorem for Besse manifolds, a generalization of the Boothby-Wang case \cite{BW}, which states that any Besse manifold is a principal $S^1$-orbibundle over a symplectic orbifold with each fiber being a Reeb orbit. Its connection 1-form is then the corresponding Besse contact form. Moreover, under a special condition, this theorem establishes a connection with Sasakian geometry. Besse manifolds have appeared in recent papers, such as \cite{Benedetti_2020} and \cite{Abbondandolo_2022}.

The idea of using symplectic homology as an invariant for contact manifolds has already been applied to a particular class of Besse manifolds known as \emph{Brieskorn manifolds}. A Brieskorn manifold is a link of a weighted homogeneous singularity, i.e., the intersection of a small sphere centered on the origin and the zero set of the corresponding polynomial, given by the following form:
$$
z_0^{a_0}+z_1^{a_1}+\cdots+z_n^{a_n}.
$$
Its link then becomes a $(2n-1)$-dimensional manifold and is denoted by $\Sigma(a_0,\dots,a_n)$. By its definition, one can obtain a canonical contact structure on the link. If one can prove that symplectic homology (or other invariants extracted from here) is independent of its filling, then it can be used as a contact invariant for manifolds in that class. For instance, Ustilovsky \cite{U1} used Brieskorn manifolds $\Sigma(p,2,2,\dots,2)$ ($4m+1$ 2's), where $p\equiv\pm 1 \pmod{8}$ and thus all diffeomorphic to the $(4m+1)$-sphere, to show that there are infinitely many contact structures on the $(4m+1)$-sphere.  Later, Fauck \cite{Fauck} used Rabinowitz-Floer homology to complete the proof, and Gutt \cite{Guttphd,Guttpesh} gave another proof using positive $S^1$-equivariant symplectic homology afterwards.

Meanwhile, a weighted homogeneous polynomial, by its definition, has a canonical $S^1$-action, and this induces an $S^1$-action on its corresponding link. If one appropriately endows a contact form for the contact structure on the link, then it will make these $S^1$-orbits into Reeb orbits. In other words, this contact manifold becomes a Besse manifold. The class of Besse manifolds is strictly larger than the class of links of weighted homogeneous singularities; for example, the Boothby-Wang bundle over a torus with Euler number $-8$ cannot be expressed as a link of hypersurface singularity (it is not even a link of complete intersection singularity \cite{simE}). We hope that our similar attempts for this strict generalization of Brieskorn cases will yield fruitful applications later.

\subsection{Setting and Main results}

We mainly focus on 3-dimensional Besse manifolds. In this case, the classification of Besse manifolds is actually reduced to (3-dimensional) Seifert fibrations, which are well-known objects in geometric topology. This observation tells us that the low-dimensional case is quite manageable. Also, since we want to have integer grading on our symplectic homology (not only for aesthetic reasons but also for needs such as the Morse-Bott spectral sequence), we shall assume that all the fillings have a trivial first Chern class. Thus, we only consider Besse manifolds with a trivial first Chern class, and we only use Stein fillings when defining our contact invariant.

Under these assumptions, we have the following theorem.

\begin{theorem}
    For a 3-dimensional Besse manifold $(M,\xi=\text{ker }\alpha)$ with a given Seifert invariants
    
    \noindent $(b,g;(\alpha_1,\beta_1),(\alpha_2,\beta_2),\cdots,(\alpha_k,\beta_k))$ whose first Chern class is trivial, the following holds:
    \begin{enumerate}[i]
        \item When the orbifold Euler characteristic is positive, the positive symplectic homology defined on its Stein filling becomes an invariant for such Besse manifolds.
        \item When the orbifold Euler characteristic is non-positive, the positive symplectic homology defined on the trivial Liouville cobordism $[0,1]\times M$ becomes an invariant for such Besse manifolds.
    \end{enumerate}
\end{theorem}
Here, the orbifold Euler characteristic is defined by,
$$
\chi:=(2-2g)-\sum^k_{j=1}\left(1-\frac{1}{\alpha_j}\right)
$$
and "the positive symplectic homology defined on the trivial Liouville cobordism" refers to the work of Cieliebak and Oancea \cite{cobo}, who constructed a variant of symplectic homology under the "hypertight condition." Due to this condition, no holomorphic planes asymptotic to a Reeb orbit exist, and thus all the Floer cylinders remain inside the cobordism, making symplectic homology well-defined there.

In fact, we will later classify all Besse manifolds with a trivial first Chern class. To establish a criterion for such structures, we construct the CW-structure for a given Besse manifold and apply obstruction theory there. Then we obtain the following criterion:

\begin{prop}\label{prev}
   The first Chern class of the Besse structure $\xi:=\ker\alpha$ of the Besse manifold $(M,\alpha)$ described by the Seifert invariant data $(0,g; (\alpha_1, \beta_1),\, \cdots \,(\alpha_{k+1},\beta_{k+1}))$ can be represented by,
   $$
     \begin{bmatrix}
        \alpha_1-1\\\alpha_2-1\\ \vdots \\\alpha_k-1\\ (2g-1)\alpha_{k+1}-1\\0
    \end{bmatrix}
    \in \Z^{k+2}/\,\,\emph{im}A_M \subseteq  H^2(M)
    $$
    where
    $$
    A_M=\begin{bmatrix}
    \alpha_1 & 0 & \cdots & 0& \beta_1 \\
0 & \alpha_2 & \cdots & 0 & \beta_2 \\
\vdots & \vdots & \ddots & \vdots & \vdots \\
0 & 0 & \cdots & \alpha_{k+1} & \beta_{k+1} \\
-1 & -1 & \cdots & -1 & 0 \\
\end{bmatrix}.
    $$
   Especially, the Besse structure has a trivial first Chern class $c_1(\xi)$ if and only if all of the following hold:
    \begin{enumerate}
        \item $\chi/e$ is an integer.
        \item For any $i$s, $\alpha_j|(\beta_j\cdot\chi/e-1)$.
    \end{enumerate}
\end{prop}

From this criterion, we can actually classify Besse manifolds with a trivial first Chern class. Such manifolds with a positive orbifold Euler characteristic (of the base orbifold) will only be the 3-sphere and links of simple singularities. When the orbifold Euler characteristic is zero, we will only have Boothby-Wang bundles over a torus with a negative Euler number.

Later, we will attempt to compute the symplectic homology for some cases, and for this, we will have to compute the Robbin-Salamon indices of periodic Reeb orbits.

\begin{theorem}[Robbin-Salamon index of a periodic Reeb orbit in a Besse manifold]
    Let $M$ be a Besse manifold of Seifert invariants $(b,g; (\alpha_1, \beta_1), \cdots  , (\alpha_k, \beta_k))$. Let $x_m$ be a principal periodic Reeb orbit of multiplicity $m$ and $x_{j,m}$ be an exceptional periodic Reeb orbit of multiplicity $m$ corresponding to $(\alpha_j,\beta_j)$ in $M$. Then the Robbin-Salamon index of each orbit is,
    \begin{align}
        &\mu(x_m)=2m\,\frac{\chi}{e}\\
        &\mu(x_{j,m})=
\begin{cases*}
    2\left\lfloor\frac{m}{\alpha_j}\cdot\frac{\chi}{e}\right\rfloor+1 & \emph{($\chi>0$)}\\
    -2\left\lfloor-\frac{m}{\alpha_j}\cdot\frac{\chi}{e}\right\rfloor-1 & \emph{($\chi<0$)}
\end{cases*}
    \end{align}
    where $\chi:=2-2g-\sum_j\left(1-\frac{1}{\alpha_j}\right)$ is the orbifold Euler characteristic and $e:=b+\sum_j \frac{\beta_j}{\alpha_j}$ is the Euler number.
\end{theorem}

\section{Preliminaries}
\subsection{Besse manifolds and Seifert fibration}
Recall that any contact manifold $M$ with a given contact form $\alpha$ admits a special vector field called the Reeb vector field.

\begin{definition}
    A \emph{Besse manifold} is a contact manifold $(M,\xi)$ that admits a contact 1-form $\alpha$ of which the Reeb flow is periodic. We call this contact structure a $\emph{Besse (contact) structure}$, and this contact 1-form a $\emph{Besse (contact 1-)form}$.
\end{definition}
\begin{remark}
    In general, all the orbits from a given flow being periodic does not mean that the flow is periodic. For example, Sullivan \cite{orbit} made a flow on a compact 5-dimensional manifold whose orbits are all periodic, but the flow itself is not. But in our case, it does due to Wadsley's theorem \cite{Wadsley}.
\end{remark}

Instead of being just an $S^1$-bundle over a symplectic manifold, any Besse manifolds are contactomorphic to an $S^1$-``orbibundle'' over a symplectic ``orbifold'' (see \cite{sasaki}, Proposition 7.1.2 and Theorem 7.1.3. and \cite{Kegel}). 

\begin{remark}
    This theorem is closely related to Sasakian geometry. A \emph{Sasakian manifold} $M$ is a contact manifold that makes the almost complex structure of SFT-type (which will be defined later) on the positive symplectization $\R_+\times M$ of $M$ integrable. The orbit space $X$ of the above theorem is K\"{a}hler if and only if $M$ is Sasakian.    
\end{remark}

Note that all of the Reeb orbits need not have the same period. Some special cases where all periods are the same are called a \emph{Zoll contact manifold} (see \cite{BW}, and also \cite{Geiges}). Any Zoll case can be seen as a $S^1$-bundle over a symplectic manifold called the \emph{Boothby-Wang Bundle} with its connection form as its contact form. Besse cases can be treated similarly. 

From now on, we only consider closed 3-dimensional cases. An $S^1$-orbibundle with a three-dimensional manifold as the total space is also called a \emph{Seifert fibration}, and this already has its classification theorem \cite{orlik}. Using this, we can give a classification theorem for 3-dimensional Besse manifolds.

\begin{theorem}[{{\cite[Cor.1.6]{Kegel}}}] \label{class}
    The classification of closed 3-dimensional Besse manifolds up to strict contactomorphisms corresponds to the classification of Seifert fibrations with closed total space and with a nontrivial (rational) Euler number (which will be explained later in detail). In other words, a closed Besse 3-manifold can be classified by the following data:
    $$
    (b,g; (\alpha_1, \beta_1),\, \cdots \, , (\alpha_k, \beta_k))
    $$
    where all are integers, $\alpha_j$ and $\beta_j$ are coprime and $0 < \beta_j \leq \alpha_j$ for all $j$s, and $e := b + \sum_j \frac{\beta_j}{\alpha_j} > 0$.
\end{theorem}

The detailed construction of each Besse manifold using these data is as follows. To construct an underlying Seifert fibration, first set a total space $\mathcal{E}$ of a trivial $S^1$-bundle $\pi:\mathcal{E}\rightarrow \Sigma$ over a genus $g$ surface $\Sigma$ with $k+1$ punctures. Each collar neighborhood near each $j$-th puncture ($j=1,\dots,k+1$) on $\Sigma$ is parametrized by polar coordinate $(r_j,\theta_j)\in(1-\varepsilon,1)\times S^1$ for some small number $\varepsilon\in(0,1)$, where its radial direction, represented by $r_j$, points outward toward the puncture. Let $p_j$s and $q_j$s be integers that satisfy
    $$
    \alpha_jq_j-\beta_jp_j=1
    $$
 and let $f_j: (0,\delta)\rightarrow(1-\varepsilon,1)$ be a smooth bijection satisfying 
    \begin{equation}\label{condi1}
    f_j'<0
    \end{equation}
where $\delta>0$ is small enough to satisfy
\begin{equation}\label{small}
1+p_j\delta^2>0 \text{ and } \beta_j-\delta^2>0
\end{equation}
for any $j$s. (Other properties that $f_j$ have shall be specified later.) Then for an open disk $D_\delta$ with radius $\delta$ and for $j=1,\dots,k+1$, attach $T_j:=D_\delta\times S^1$ to each collar neighborhood $(1-\varepsilon,1)\times S^1\times S^1$ of each boundary component of $\mathcal{E}$ via the following map:
    \begin{align}
    T_j\supseteq(0,\delta)\times S^1\times S^1&\rightarrow(1-\varepsilon,1)\times S^1\times S^1\subseteq\Sigma\times S^1 \nonumber  \\
    (r_j,\theta_j,\varphi_j)&\mapsto(s_j=f_j(r_j),u_j=-\alpha_j\theta_j-p_j\varphi_j,t=\beta_j\theta_j+q_j\varphi_j).\label{gluingmap}
\end{align}
Here, $(\alpha_{k+1},\beta_{k+1})$ is defined by $(1,b)$. 

\begin{remark}\label{redundancy}
When we defined integers $p_j$s and $q_j$s, it only required to satisfy $
    \alpha_jq_j-\beta_jp_j=1
    $
and thus instead of a given pair $(p_j,q_j)$, one can replace it with a pair $(p_j,q_j)+\epsilon(\alpha_j,\beta_j)$ for any integer $\epsilon$. This, in our construction, is related to the freedom of parametrizing each $T_j$s. Since there are the following diffeomorphisms
    \begin{align*}
    D_\delta\times S^1&\rightarrow D_\delta\times S^1\\
    (r,e^{i\theta},e^{i\varphi})&\mapsto ((r,e^{i(\theta+\epsilon\varphi)},e^{i\varphi})),
    \end{align*}
one can have different parametrizations of $T_j$ via those maps, and once one chooses a pair $(p_j,q_j)$ and a parametrization, the rest of the possible maps exactly correspond to each possible pair of $(p_j,q_j)$s respectively (up to orientation of $T_j$ or sign change of $\epsilon$.)

In fact, when looking at the whole Besse manifold $M$, one can ignore this redundancy of choice of a pair $(p_j,q_j)$. This is because if one interprets the attachment $T_j$s attaching a pair of a 2-handle and a 3-handle, one only needs $\alpha_j$ and $\beta_j$ to specify the attaching map for the 2-handle. Since any data on the attachment of the 3-handle gives the same smooth manifold, one does not need $p_j$ and $q_j$ to examine the total handlebody. 
    \end{remark}

Since we now have our underlying Seifert fibration, we also describe how to endow the actual Besse contact structure on this manifold. For this, we choose contact forms on each part in the above construction that gives the Reeb flow we want, and show that these choices are compatible with each other and give a (global) Besse 1-form for the whole space. 
For the Besse contact form on $\mathcal{E}$, we need the following lemma.
\begin{lemma}
    For any compact Riemann surface $\Sigma$ with $k+1$ boundary components and for any real numbers $a_j(j=1,\dots,k+1)$ of which the total sum is positive, there exists a 1-form $\lambda$ such that $d\lambda$ is an area form of $\Sigma$, its restriction near the $j$-th $(j=1,\dots,k+1)$ boundary component is 
    $$
    \rho_j(s_j)du_j.
    $$
    Here, $(s_j,u_j)$ is the coordinate system we have defined above on $\Sigma$ near the $j$-th puncture for all $j=1,\dots,k+1$, where $\{s_j=1\}$ is the boundary component, and functions $\rho_j:(1-\varepsilon,1)\rightarrow \R$ $(j=1,\dots,k+1)$ are given smooth functions satisfying $\rho_j'>0$ and $\lim_{s_j\nearrow{1}}\rho_j(s_j)=a_j$.
\end{lemma}
\begin{proof}
Modify the initial part of the Thurston-Winkelnkemper construction \cite{ThurstonWinkelnkemper}. Since we assumed that the total sum of $a_j$s is positive, the construction works similarly well.
\end{proof}
Now, choose our functions $\rho_j(s_j)$ that satisfy
\begin{equation}\label{condi2}
    \rho_j'>0
\end{equation}
and $\lim_{s_j\nearrow{1}}\rho_j(s_j)=\beta_j/\alpha_j$. Due to our assumption $e:=b+\Sigma_j\beta_j/\alpha_j>0$, we can use the above lemma for these functions and get a 1-form $\lambda$ on $\Sigma$ from the lemma. Then choose a contact form on $\mathcal{E}$ by
$$\lambda+dt.$$
Here, $\lambda$ is the abused notation for the pullback 1-form $\pi^*\lambda$. Then one can assume that the restriction of the above contact form to a collar neighborhood $(1-\varepsilon,1)\times S^1\times S^1$ of the $j$-th boundary component of $\Sigma\times S^1$ becomes,
$$\rho_j(s_j)du_j+dt.$$
On the $j$-th solid torus $T_j$, we choose 
$$
r_j^2d\theta_j+\frac{1}{\alpha_j}(1+p_jr_j^2)d\varphi_j
$$
as a contact form. Note that not only the 1-forms we choose as contact forms are really contact forms, it also gives the Reeb flow we want on each part. The Reeb vector field on $\Sigma\times S^1$ is $\frac{\partial}{\partial t}$ and on the $j$-th solid torus is 
$$
-p_j\frac{\partial}{\partial \theta_j}+\alpha_j\frac{\partial}{\partial \varphi_j}.
$$
The only thing left is to adjust the function $f_j(r_j)$ to associate these Besse 1-forms into a single Besse 1-form on the Besse manifold $M$. For this, one only needs to assume
\begin{equation}\label{condi3}
 \rho_j\,\circ f_j(r_j)=\frac{1}{\alpha_j}(\beta_j-r_j^2)
\end{equation}
for the functions $f_j(r_j),\rho_j(s_j)$, and one can actually find such a function $f_j$ satisfying conditions (\ref{condi1}) and (\ref{condi3}) because $\rho_j$ strictly increases ((\ref{condi2})).

\begin{figure}
    \centering
    \includegraphics[width=1\linewidth]{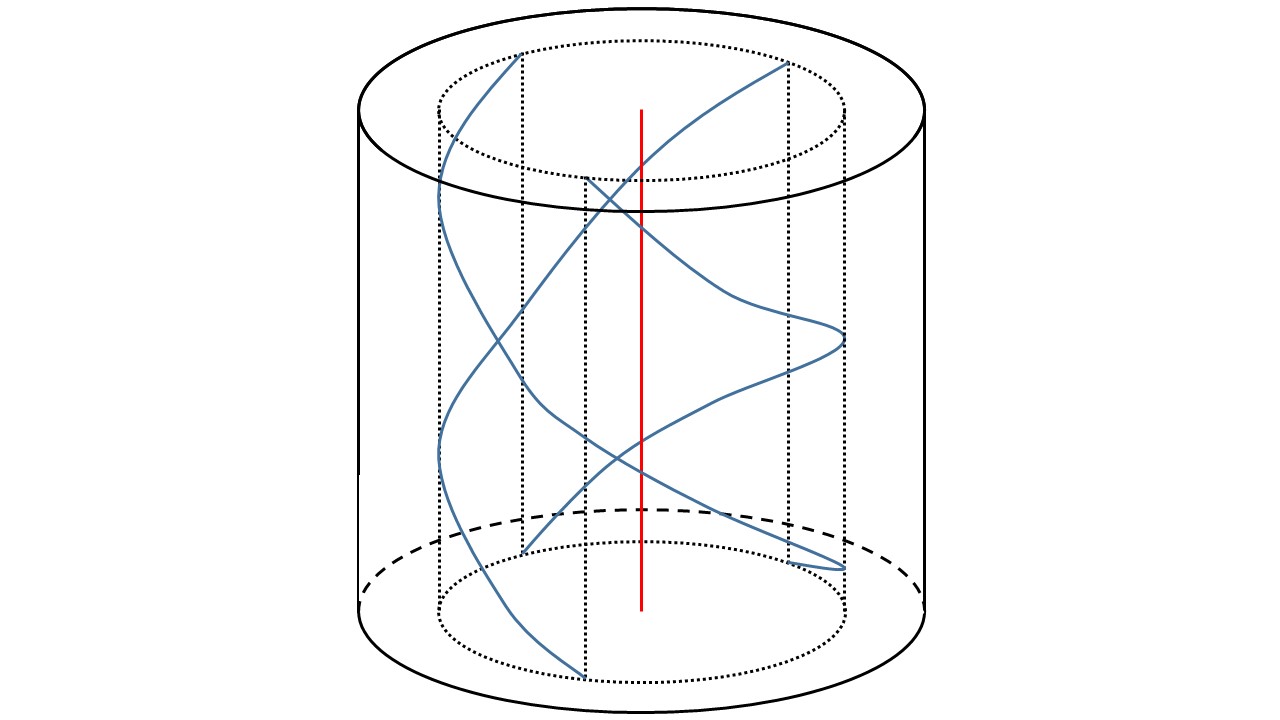}
    \caption{A tubular neighborhood of an exceptional orbit of (3,2). (The top face and the bottom face are identified.) The red line is the exceptional orbit, and the blue lines are some of the principal orbits nearby.}
    \label{except}
\end{figure}

\begin{remark}
    For each solid torus $T_j$, the Reeb flow can be described as the following $S^1$-action:
\begin{align}\label{action}
\begin{split}
     S^1\times D^2\times S^1&\rightarrow D^2\times S^1 \\
    (u,x,v)\,&\mapsto (u^{-p_j}\cdot x,u^{\alpha_j}\cdot v).
\end{split}
\end{align}
The orbit $\{x=0\}$ in this torus (the red orbit in Figure \ref{except}) is called the \emph{exceptional orbit} of $(\alpha_j,\beta_j)$ and orbits besides that are called \emph{principal orbits}. Then one can see from the above action that near the exceptional orbit, the principal orbits wrap around it in the form of $(\alpha_j,-p_j)$-torus knots. 
\end{remark}

\begin{remark}\label{free}(\cite{orlik})
    Although we restricted $\beta_j$s to certain conditions, one can instead set $\beta_j$s to be just an integer. Then the following three sets of data give the same Besse manifold:
    $$
    (b,g; (\alpha_1, \beta_1),\, \cdots \,(\alpha_j,\beta_j)\, \cdots \, , (\alpha_k, \beta_k))
    $$
    \begin{center}
        and
    \end{center}
    $$(b+1,g; (\alpha_1, \beta_1),\, \cdots \,(\alpha_j,\beta_j-\alpha_j)\, \cdots \, , (\alpha_k, \beta_k)).
    $$
Especially, one can make $b$ to be 0 and change some other $\beta_j$s. In addition, inserting/deleting a pair $(1,0)$ will give the same result. 
\end{remark}

\begin{remark}
    Unlike other invariants, the Seifert invariant $b$ itself is difficult to detect geometrically. This is because a surgery along a fiber in a trivial $S^1$-bundle over a closed surface, according to the pair $(1,b)$ will again be a trivial $S^1$-bundle. But the Euler number $e$, the total sum of $b$ and all $\beta_j/\alpha_j$s, which we will talk about later, will be detected in the form of topological invariants like, for instance, homology/cohomology groups. 
\end{remark}


Now, in this particular setting, let us introduce some useful invariants. One would be the \emph{Euler number} $e$ of a Seifert fibration $M\rightarrow S$, which was aforementioned in Theorem \ref{class}. $$e(M):=b+\sum_j \frac{\beta_j}{\alpha_j}$$ This number is basically the generalization of the Euler number for genuine $S^1$-bundles: Both coincide for $S^1$-bundle cases, and in fact, the prior has the same role as the original Euler number did. For example, the Euler number of a Seifert fibration can be seen as a ``self-intersection number of the zero section of the disc bundle induced from the Seifert fibration'' \cite{seifertt}.

\begin{remark}\label{ori}
    The condition $e>0$ used in Theorem \ref{class} is related to the orientation of $M$. A Seifert fibration from the data $(b,g;(\alpha_1,\beta_1),(\alpha_2,\beta_2),\cdots,(\alpha_k,\beta_k))$ is orientation-reversing diffeomorphic to $(-k-b,g;(\alpha_1,\alpha_1-\beta_1),(\alpha_2,\alpha_2-\beta_2),\cdots,(\alpha_k,\alpha_k-\beta_k))$. Note that the Euler numbers of each Seifert fibrations have only a sign difference; thus, only one of them with a positive $e$ can be a Besse manifold. In our construction of the Besse manifold, this appears as a sufficient condition for the existence of the 1-form $\lambda$ we used.  
    
    Meanwhile, in \cite{Kegel}, the authors seemed to allow contact structures with opposite orientation on the same Seifert fibration, stating that $e$ being nonzero is enough.
\end{remark}

Another invariant to introduce would be the \emph{orbifold Euler characteristic}, which is a rational number for the base orbifold $X$ of a Seifert fibration $M\rightarrow X$. For a fibration from the data $(b,g; (\alpha_1, \beta_1),\, \cdots \, , (\alpha_k, \beta_k))$, the orbifold Euler characteristic is defined by \cite{scott}:
$$
\chi(S):=2-2g-\sum_j\left(1-\frac{1}{\alpha_j}\right).
$$
This invariant generalizes the original Euler characteristic in the following way: If a 2-dimensional closed orbifold has a $d$-fold covering orbifold, the orbifold Euler characteristic of the covering equals $d$ times the orbifold Euler characteristic of the other \cite{scott}.

From the perspective of symplectic geometry, there is the following Stein-fillability result.
\begin{prop}\label{BesseStein}
    Any 3-dimensional Besse manifold is Stein fillable.
\end{prop}
\begin{proof}
    Let $(M,\alpha,\xi:=\ker\alpha)$ be a Besse manifold. Since $(\xi,d\alpha)$ is a rank 2 symplectic vector bundle, we have a complex structure $J$ on $\xi$. We also get a metric $g(\cdot,\cdot):=d\alpha(\cdot,J\cdot)$. By averaging this $g$ along the $S^1$-action generated by the Reeb vector field, we can have an $S^1$-invariant metric $\Bar{g}$ and get the complex structure $\Bar{J}$ such that $\Bar{g}(\cdot,\cdot):=d\alpha(\cdot,\Bar{J}\cdot)$. Since the Reeb vector field is transverse to $\xi$, by Theorem 2.1 of \cite{CR} we get a pseudoconvex holomorphic filling. Blow this filling down into minimal filling and apply Theorem 5.64 of \cite{CE} to get a Stein filling.
\end{proof}

\subsection{Examples of Besse manifolds}

\begin{example} [$S^3$]
The standard contact structure over $S^3:=\{(z_1,z_2)\in\C^2:|z_1|^2+|z_2|^2=1\}$ can be described as a kernel of the 1-form 
$$
\alpha_{std}:=\left. \frac{1}{2}(-y_1dx_1+x_1dy_1-y_2dx_2+x_2dy_2)\right|_{S^3}
$$
where $z_j=x_j+iy_j$ $(j=1,2)$. The Reeb vector field of this 1-form gives exactly the $S^1$-action induced from the one in $\C^2$, described as
$$
e^{it}\cdot(z_1,z_2):=(e^{it}z_1,e^{it}z_2)\,\,\,(t\in\R).
$$
One can easily see that all the orbits of this action are periodic with the same period $2\pi$. (Note that this action gives the Hopf fibration, the $S^1$-bundle over the sphere.)

One can get the Seifert invariants from direct computation. Let $T_j\subseteq S^3$ $(j=1,2)$ be 
$$
\mathcal{T}_j:=\left\{(z_1,z_2)\in S^3:|z_j|^2\leq \frac{1}{2}\right\}.
$$
Once we trivialize $\mathcal{T}_1$ and $\mathcal{T}_2$ as trivial $S^1$-bundles over a disc of radius $(1/2)^{1/2}$ by 
\begin{align*}
    (x,e^{i\varphi})&\mapsto(xe^{i\varphi},(1-|x|^2)^{1/2}e^{i\varphi})\in \mathcal{T}_1\\
    (x,e^{i\varphi})&\mapsto((1-|x|^2)^{1/2}e^{i\varphi},xe^{i\varphi})\in \mathcal{T}_2,
\end{align*}
the transition map on $T:=\partial \mathcal{T}_1$ from $\mathcal{T}_1$ to $\mathcal{T}_2$ becomes
$$
(e^{i\theta},e^{i\varphi})\mapsto (e^{-i\theta},e^{i(\theta+\varphi)})
$$
where $\theta$ is the angular coordinate of the disc. Note that the $\theta$-coordinate on $\mathcal{T}_2$ is opposite to what we used during the construction of Seifert fibrations. Adjusting this will give the Seifert invariant as $b=1$, $g=0$, and no exceptional fibers.

Another way to give a Besse contact structure on $S^3$ would be the \emph{ellipsoid}. The underlying manifold is defined by
$$
E_{a,b}:=\left\{(z_1,z_2)\in\C^2:\frac{|z_1|^2}{a}+\frac{|z_2|^2}{b}=1\right\}
$$
where again $z_j=x_j+iy_j$, and $a,b$ are coprime integers. (Obviously, this is diffeomorphic to $S^3$.) The contact form on $E_{a,b}$ is defined by
$$
\alpha_{a,b}:=\left. \left(\frac{1}{2}(-y_1dx_1+x_1dy_1)+\frac{1}{2}(-y_2dx_2+x_2dy_2)\right)\right|_{E_{a,b}}.
$$
The $S^1$-action given by the Reeb vector field becomes,
$$
e^{it}\cdot(z_1,z_2):=(e^{it/a}z_1,e^{it/b}z_2)\,\,\,(t\in\R).
$$
All the orbits for this action are again periodic with period $2\pi\cdot ab$, except for two exceptional orbits $\{z_1=0\}$ and $\{z_2=0\}$ of which each period is $2\pi b$ and $2\pi a$, respectively. The Seifert invariants for this one become,
$$
\left(\frac{1-a^*a-b^*b}{ab},0;(a,b^*),(b,a^*)\right)
$$
where $a^*$ (resp. $b^*$) is an integer satisfying $a^*a\equiv1$ (mod $b$) (resp. $b^*b\equiv1$ (mod $a$)).
\end{example}

\begin{example}[Link of weighted homogeneous singularity \& simple singularities] \label{WH}
A polynomial $f(z_0,z_1,\cdots,z_n)\in\C[z_0,z_1,\cdots,z_n]$ is called \emph{weighted homogeneous polynomial} when there exist positive integers $w_0,w_1,\cdots,w_n,N$ such that for any $t\in\C$,
$$
f(t^{w_0}z_0,t^{w_1}z_1,\cdots,t^{w_n}z_n)=t^Nf(z_0,z_1,\cdots,z_n)
$$
holds. Note that $f(0)=0$. Now, assume $n=2$ and let $V$ be the variety of this polynomial. We call the germ of $V$ near $0$ \emph{the weighted homogeneous singularity}. Then we define \emph{the link of the singularity} as the intersection of a ball with a small radius $\varepsilon>0$ centered on the zero and the variety $V$. In other words, it's the set of points $M$ (which is a manifold) satisfying the following equations:

\begin{align*}
        f(z_0,z_1,z_2)=0&\\
    |z_0|^2+|z_1|^2+|z_2|^2&=\varepsilon^2.
\end{align*}
$M$ then becomes a contact manifold naturally by taking $TM\cap JTM$ as its contact structure, where $J$ is the (almost) complex structure on the variety $V$. Also, from these equations defining $M$, there is a natural $S^1$-action defined by
$$
e^{it}\cdot(z_0,z_1,z_2):=(e^{iw_0t}z_0,e^{iw_1t}z_1,e^{iw_2t}z_2)
$$
and one can actually show that this same $S^1$-action can be induced from the Reeb flow by restricting the following 1-form 
$$
\frac{1}{\varepsilon^2w_0} (-y_0dx_0+x_0dy_0)+\frac{1}{\varepsilon^2w_1} (-y_1dx_1+x_1dy_1)+\frac{1}{\varepsilon^2w_2} (-y_2dx_2+x_2dy_2)
$$
on $\C^3$ to $M$ and take that restricted 1-form as the contact form. This contact form will then endow a Besse structure on $M$. The period of each principal orbit will be $2\pi$, but orbits like ones in $\{z_0=0\},\{z_1=0\},\{z_2=0\}$ may have different periods. The periods of each type will be $2\pi/g_{12},2\pi/g_{02}$, and $2\pi/g_{01}$, respectively where $g_{ij}$ means $gcd(w_j,w_j)$.

There are certain types of weighted homogeneous singularities named \emph{simple singularities} (or \emph{ADE-type singularities}) given by the following polynomials:
\begin{align*}
    A_l & :z_0^2+z_1^2+z_2^{l+1} \\
    D_l & :z_0^2+z_1^2z_2+z_2^{l-1} \,\,\, (n\geq4)\\
    E_6 & :z_0^2+z_1^3+z_2^4 \\
    E_7 & :z_0^2+z_1^3+z_1z_2^3 \\
    E_8 & :z_0^2+z_1^3+z_2^5
\end{align*}
One can compute the Seifert invariants \emph{\cite{orlik}} as in the following table.

\begin{table}[h!]
\centering
 \begin{tabular}{|c | c | c | c | c |} 
 \hline
 Type & $(w_0,w_1,w_2,N)$ & $b$ & $g$ & $(\alpha_j,\beta_j)$ \\ 
 \hline\hline
 $A_l$ ($l$:odd) & $(\frac{l+1}{2},\frac{l+1}{2},1,l+1)$ & -2 & 0 & $(\frac{l+1}{2},\frac{l-1}{2}),(\frac{l+1}{2},\frac{l-1}{2})$\\ 
 $A_l$ ($l$:even) & $(l+1,l+1,2,2(l+1))$ & -1 & 0 & $(l+1,\frac{l}{2}),(l+1,\frac{l}{2})$\\ 
 $D_l$ & $(l-1,l-2,2,2(l-1))$ & -1 & 0 & $(2,1),(2,1),(l-2,1)$\\ 
 $E_6$ & $(6,4,3,12)$ & -1 & 0 & $(2,1),(3,1),(3,1)$\\ 
 $E_7$ & $(9,6,4,18)$ & -1 & 0 & $(2,1),(3,1),(4,1)$\\ 
 $E_8$ & $(15,10,6,30)$ & -1 & 0 & $(2,1),(3,1),(5,1)$\\ 
 \hline
 \end{tabular}
 \caption{Weights and Seifert invariants for links of simple singularities}
\end{table}
\end{example}

\subsection{Symplectic Homology}

Let $(W,M:=\partial W,\omega=d\lambda)$ be a $2n$-dimensional Liouville domain. ($(M,\alpha:=\lambda|_{M})$ is a contact manifold.) Symplectic homology \cite{SH1} \cite{SH2} is an invariant defined on a Liouville domain as ``a Morse homology of the (contractible) free loop space over $W$." For more details about the construction of this invariant, one can look for references like \cite{AD}, \cite{Abou}, and \cite{Ritter}. We also use other Floer-type theories, like positive symplectic homology or symplectic homology on Liouville cobordism \cite{cobo}.Our convention for symplectic homology and positive symplectic homology basically follows that of \cite{bourgeois2009symplectic}:
\begin{itemize}
    \item The Liouville domain is completed by attaching the positive symplectization of the contact boundary $M$. Let this completion $W$. The positive symplectic part will often be called ``the infinity cone":$(\R_{\geq1}\times\partial W,d(r\alpha)$
    \item The Hamiltonian used in our setting is $\mathcal{C}^2$-small in $W$ and is an increasing function that depends only on the radial coordinate. This Hamiltonian is linear at the infinity cone with slope $s(H)$, which is not a period of a certain periodic Reeb orbit. We will call these Hamiltonian ``admissible."
    \item The Hamiltonian vector field of a Hamiltonian $H$ is defined by $\omega(X_H,\cdot)=dH$.
    \item Our action functional is $\mathcal{A}_H(x):=-\int_{S^1} x^*\lambda-\int^1_0 H(x(t),t)dt$.
    \item All almost complex structures used here are of SFT type, i.e. 
    \begin{align*}
        J&(r\frac{\partial}{\partial r})=R\\
        J&R=-r\frac{\partial}{\partial r}\\
        J&\xi=\xi.
    \end{align*}
    Thus, on the infinity cone, where the Hamiltonian only depends on the radial coordinate (say $H(r,x)=:h(r)$), the Hamiltonian vector field becomes $$X_H=-h'(r)R.$$
    \item The Floer equation is defined as follows:
    \begin{align}
    \begin{split}\label{Fl}
        \partial_su+&J(\partial_tu-X_H)=0 \\
        \lim_{s\rightarrow \pm \infty}& u(s,\cdot) = x_{\pm}.
    \end{split}
\end{align}
    \item For any generator $x$ of the Floer chain complex, the grading of $x$ is defined as $|x|:=-\mu_{CZ}(x)$.
    \item The differential $\partial:CF_*(W,H;\mathbb{Z}/2\mathbb{Z})\rightarrow CF_*(W,H;\mathbb{Z}/2\mathbb{Z})$ is defined as follows:
$$
\partial x_-:=\sum_{\substack{x_+\in \mathcal{P}(H)\\|x_-|-|x_+|=1}}|\mathcal{M}(x_-,x_+;H,J)/\mathbb{R}|\cdot x_+ (\text{mod}2),
$$
where $\mathcal{M}(x_-,x_+;H,J)$ is the moduli space of Floer cylinders $u$ with finite energy, i.e. 
$$
\int_{\R\times S^1}\omega(\partial_su,J\partial_su)<\infty
$$
asymptotically converges to orbits $x_\pm$;
$$
\lim_{s\rightarrow\pm\infty}u(s,\cdot)=x_\pm.
$$

\item For the continuation map from $(H_-,J_-)$ to $(H_+,J_+)$ where $s(H_-)\leq s(H_+)$, construct a monotone increasing homotopy $H_s$ of admissible Hamiltonians between $H_-$ and $H_+$ and a homotopy of almost complex structure $J_s$ of SFT type between $J_-$ and $J_+$. Then define a moduli space $\mathcal{M}_s(x_-,x_+;H,J)$ of solutions of s-dependent Floer equation
\begin{align}
    \begin{split}
        \partial_su+&J_s(\partial_tu-X_{H_s})=0 \\
        \lim_{s\rightarrow \pm \infty}& u(s,\cdot) = x_{\pm}.
    \end{split}
\end{align}
Here, $x_\pm$ is a Hamiltonian orbit for $H_\pm$, respectively.
\item A continuation map $\phi_{H_-,H_+}:CF_*(W,H_-;\mathbb{Z}/2\mathbb{Z})\rightarrow CF_*(W,H_+;\mathbb{Z}/2\mathbb{Z})$ from $(H_-,J_-)$ to $(H_+,J_+)$ is given by,
$$
\phi_{H_-,H_+} (x_-):=\sum_{\substack{x_+\in \mathcal{P}(H)\\|x_+|=|x_-|}}|\mathcal{M}_s(x_-,x_+;H,J)|\cdot x_+ (\text{mod}2)
$$
\end{itemize}

\section{Computation of $c_1$ of a Besse structure}
In this section, we will provide a cell (or handle) decomposition of the total space $M$ of a Seifert fibration constructed from the data $(b,g; (\alpha_1, \beta_1),\dots,(\alpha_k,\beta_k))$ in order to compute the first Chern class of the Besse structure. Instead of using the Seifert invariant $b$, we will occasionally assume $b=0$, considering our data as $(0,g; (\alpha_1, \beta_1),\dots,(\alpha_{k+1},\beta_{k+1}))$ without imposing the condition $0\leq\beta_j\leq\alpha_j$. (See Remark \ref{free}.)

\begin{figure}
    \centering
    \includegraphics[width=1\linewidth]{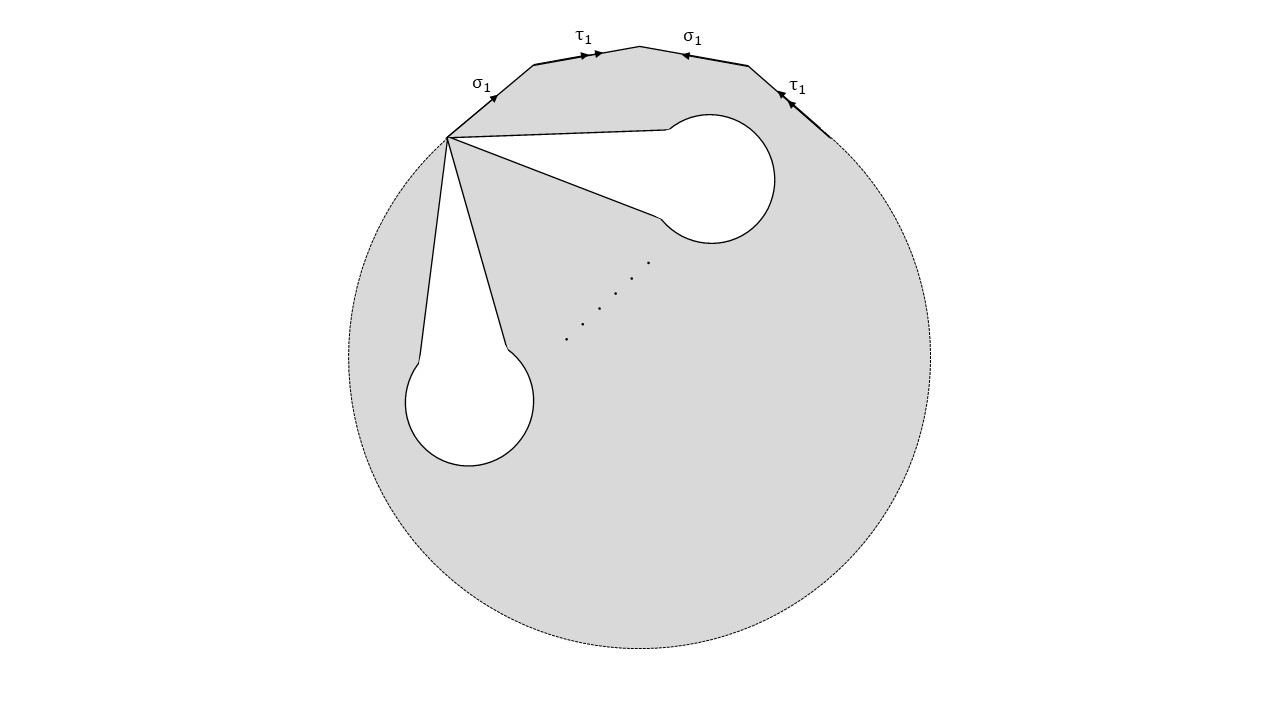}
    \caption{CW-complex structure of $\Sigma$. \\Each pair of edges of $4g$-gon with $k+1$ punctures is identified along the marked direction.}\label{surf}
\end{figure}

Recall that a Seifert fibration is obtained by performing multiple surgeries on each fiber of the trivial $S^1$ bundle over a surface of genus $g$ according to our data. Let $\Sigma$ be a compact Riemann surface of genus $g$ with $k+1$ boundary components. This surface with boundary can be represented as a CW-complex by taking a 0-cell, $2g+k+1$ 1-cells (where $2g$ of them come from the genus and $k+1$ from the boundary components), and a 2-cell. (Figure \ref{surf})
Then, by taking the product cell structure, the total space $\mathcal{E}$ of the trivial $S^1$-bundle $\pi:\mathcal{E}\rightarrow \Sigma$ over $\Sigma$ as a CW-complex consists of a 0-cell, $2g+k+2$ 1-cells, $2g+k+2$ 2-cells, and a 3-cell. Label the 0-cell and the 3-cell as $e_0$ and $e_3$. Among the $2g+k+2$ 1-cells, $2g$ of them, labeled $\sigma_l,\tau_l$ $(l=1,\dots,g)$, are from the genus, $k+1$ of them, labeled $\mu_j$ $(j=1,\dots,k+1)$, are from the boundary components of $\Sigma$, and the remaining one, labeled $\lambda$, is from a fiber $S^1$. Among the $2g+k+2$ 2-cells, $2g$ of them, labeled $\mathfrak{S}_l,\mathfrak{T}_l$ $(l=1,\dots,g)$, are from products of 1-cells in $\Sigma$ from the genus and $S^1$, $k+1$ of them, labeled $\mathfrak{M}_j$ $(j=1,\dots,k+1)$, are from products of 1-cells in $\Sigma$ from its boundary components and $S^1$ (hence each of them will be the boundary components of $\mathcal{E}$), and the remaining one, labeled $\mathfrak{B}$, is from the 2-cell of $\Sigma$.

On this CW-complex, attach $k+1$ 2-cells $C_j$s ($j=1,2,\dots, k+1$) along the torus knot on the torus boundary component that wraps along the meridian direction $\alpha_j$ times and the longitudinal direction $\beta_j$ times. Then attach $k+1$ 3-cells $D_j$s ($j=1,2,\dots, k+1$) on each remaining boundary component to obtain our Seifert fibration.

From this CW-complex structure, one can compute the homology/cohomology of this Besse manifold. We focus on computing the second homology/cohomology group. First, we set our ordered basis for the cellular complex as follows: $$\{\mu_1,\dots,\mu_{k+1},\lambda,\sigma_1,\tau_1,\dots,\sigma_g,\tau_g\}$$ for 1-cells, $$\{C_1,\dots,C_{k+1},\mathfrak{B},\mathfrak{M}_1,\dots,\mathfrak{M}_{k+1},\mathfrak{S}_1,\mathfrak{T}_1,\dots,\mathfrak{S}_g,\mathfrak{T}_g\}$$ for 2-cells, and $$\{D_1,\dots,D_{k+1},e_3\}$$ for 3-cells. Then the boundary map from the second degree to the first is
\[
\begin{bmatrix}
\begin{matrix}
\alpha_1 & 0 & \cdots & 0& -1 \\
0 & \alpha_2 & \cdots & 0 & -1 \\
\vdots & \vdots & \ddots & \vdots & \vdots \\
0 & 0 & \cdots & \alpha_{k+1} & -1 \\
\beta_1 & \beta_2 & \cdots & \beta_{k+1} & 0 \\
\end{matrix} & \rvline & \bigzero \\
\hline
\bigzero & \rvline & \bigzero
\end{bmatrix},
\]
(Note that -1's on the $(k+1)$-th column are due to our choice of orientations of $\mathfrak{B}$ and 1-cells, so that $\mathfrak{B}$ has boundary components homologous to $-\mu_j$s.) and the boundary map from the third degree to the second is
\[
\begin{bmatrix}
\bigzero \\
\hline
\begin{matrix}
1 & 0 & \cdots & 0& -1 \\
0 & 1 & \cdots & 0 & -1 \\
\vdots & \vdots & \ddots & \vdots & \vdots \\
0 & 0 & \cdots & 1 & -1 \\
\end{matrix} \\
\hline
\bigzero
\end{bmatrix},
\]
(Similarly as above, note that -1's on the $(k+2)$-th column are due to our choice of orientation relation between $\mathfrak{B}$ and 2-cells.) where the first $k+2$ rows are zero. Thus, using the boundary map we previously wrote, the degree 2 cohomology group of our Seifert fibration is
\[
H^2(M)\simeq \mathbb{Z}^{2g}\oplus\left(\mathbb{Z}^{k+2}/\text{im}A_M\right)
\]
where $A_M$ is the matrix
\begin{equation}
    A_M=\begin{bmatrix}
    \alpha_1 & 0 & \cdots & 0& \beta_1 \\
0 & \alpha_2 & \cdots & 0 & \beta_2 \\
\vdots & \vdots & \ddots & \vdots & \vdots \\
0 & 0 & \cdots & \alpha_{k+1} & \beta_{k+1} \\
-1 & -1 & \cdots & -1 & 0 \\
\end{bmatrix}.
\end{equation}

Note that the determinant of $A_M$ equals $\alpha_1\alpha_2\cdots\alpha_k e$ (where $e$ was the Euler number of the Seifert fibration defined at Theorem \ref{class}), so for any given Besse manifold $M$, the latter summand of $H^2(M)$ is its torsion part. The generators of this torsion part are exactly the dual elements corresponding to $C_j$s.

One can now discuss the first Chern class of our Besse structure $(M,\alpha)$. For this, we first symplectically trivialize $\xi|_{\mathcal{E}}$ and $\xi|_{T_j}$s. (From now on, every trivialization refers to ``symplectic" trivialization.) Recall that the Besse form on $\mathcal{E}$ is $\lambda+dt$, where $\lambda$ is a 1-form satisfying some specific conditions and $t$ is the coordinate for the fiber-wise direction. If we let $\Phi_\Sigma:T\Sigma\rightarrow \Sigma\times\R^2$ be the trivialization of the tangent bundle $T\Sigma$ over $\Sigma$, then one can trivialize $\xi|_{\mathcal{E}}$ 
by,
\begin{align*}
    \Phi:\xi|_{\mathcal{E}}&\rightarrow \mathcal{E}\times \R^2\\
    (p,e^{it},v)&\mapsto(p,e^{it},\Phi_\Sigma(\pi_*v)).
\end{align*}
In fact, the symplectic form of the Besse structure $\xi$ is $d\lambda$, the pullback of an area form of $\Sigma$, so the above trivialization is symplectic. 

\begin{figure}
    \centering
    \includegraphics[width=1\linewidth]{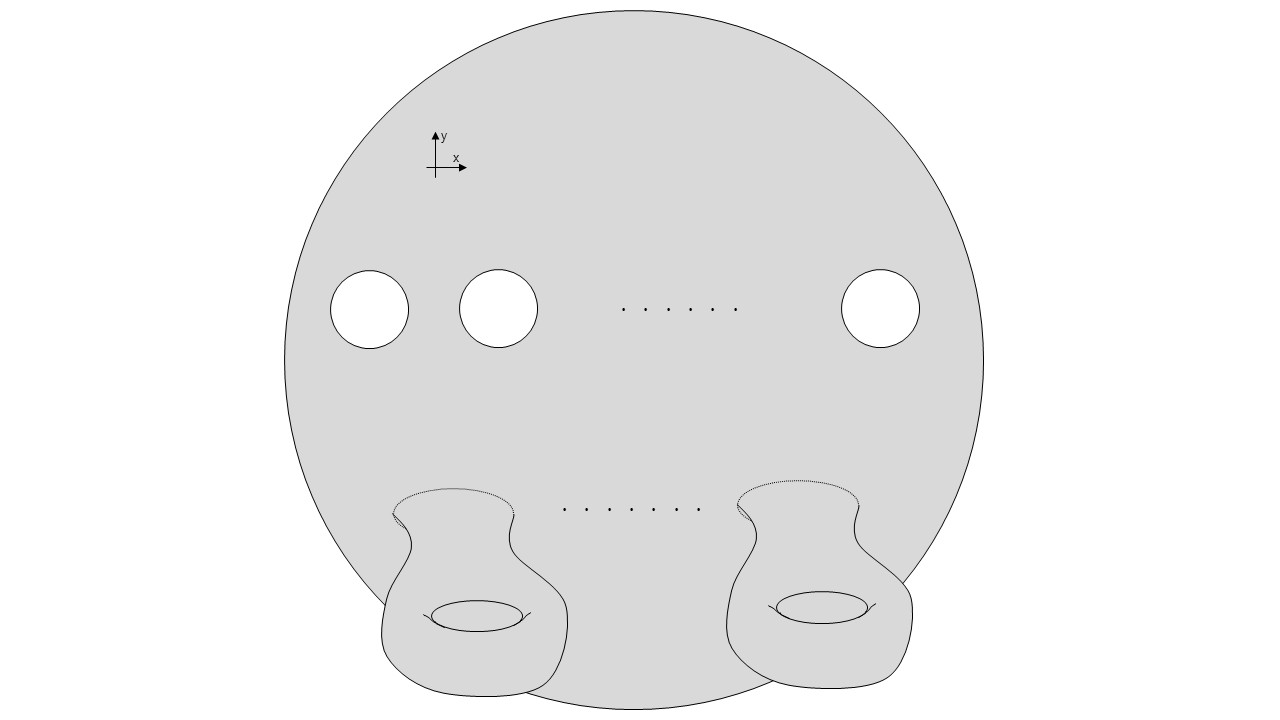}
    \caption{$\Sigma$ as a union of the punctured disc $D$ (with its Cartesian coordinate) and punctured tori $\Sigma_l$s.}
    \label{SigmaTriv}
\end{figure}

\begin{figure}
    \centering
    \includegraphics[width=1\linewidth]{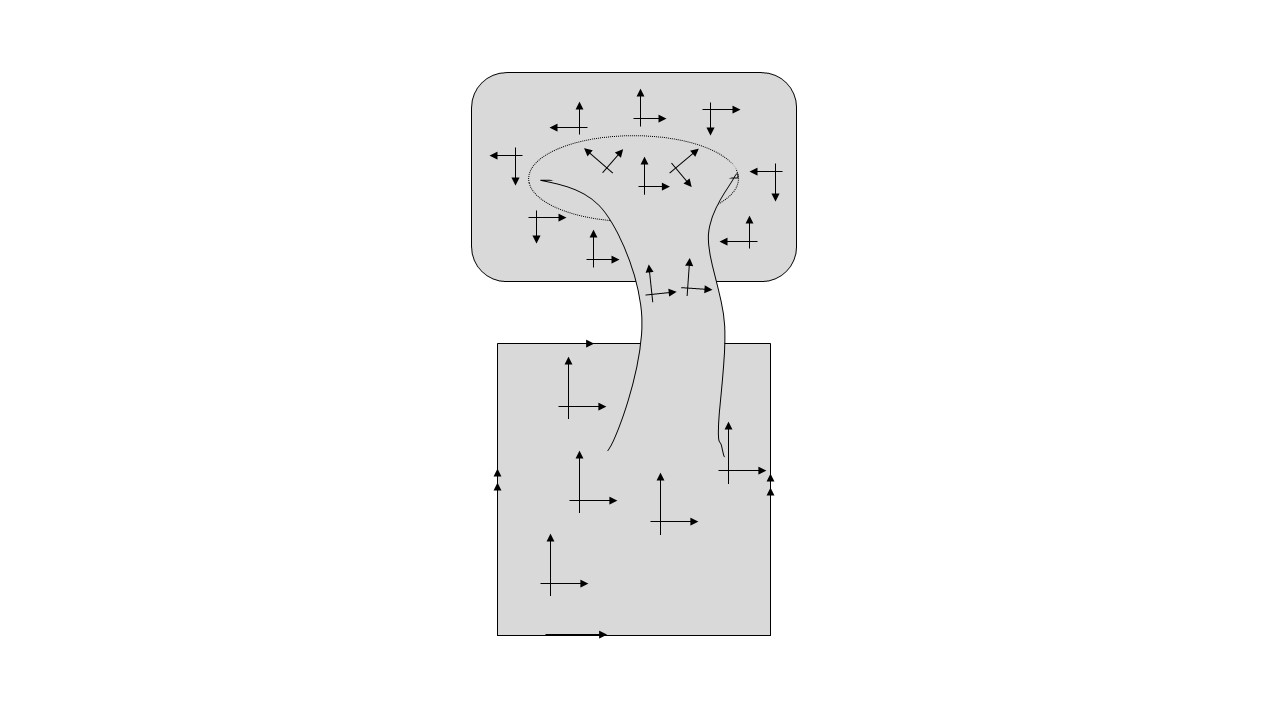}
    \caption{Trivialization on a torus with boundary and the transition map from it.}
    \label{ASingleTorus}
\end{figure}

So, to trivialize $\xi|_{\mathcal{E}}$, we only need to trivialize $T\Sigma$. First, split $\Sigma$ into a disc $D$ with $k+g$ punctures and $g$ tori $\Sigma_l$s ($l=1,\dots,g$) with each having a connected boundary. (Figure \ref{SigmaTriv}) Then trivialize the bundle restricted on the punctured disc
by identifying the disc as a subset of $\R^2$ and using its Cartesian coordinate. Let this trivialization 
\begin{align*}
\Phi_D:TD&\rightarrow D\times\R^2\subseteq T\R^2\\
\frac{1}{\Lambda(p)}\cdot\frac{\partial}{\partial x}_p,\frac{\partial}{\partial y}_p&\mapsto (p,(1,0)),(p,(0,1)).
\end{align*}
Here, 
$$
\Lambda(p):=d\lambda\left(\frac{\partial}{\partial x},\frac{\partial}{\partial y}\right).
$$
If we temporarily identify each torus by rectangles with each pair of non-adjacent edges identified, we can also trivialize
$T\Sigma_l$ by identifying this rectangle as a subdomain of $\R^2$ and induce trivialization from the canonical one of $\R^2$. (Figure \ref{ASingleTorus}) Let this trivialization $\Phi_{\Sigma_l}:T\Sigma_l\rightarrow\Sigma_l\times\R^2$. 
Parameterize each collar neighborhood of the boundary components on $D$ and $\Sigma_l$ in polar coordinates by annulus $A(1,1+\varepsilon)$ $(\varepsilon\in(0,1))$. The gluing of $\Sigma_l$ on $D$ then becomes,\begin{align*}
    F_l:A(1,1+\varepsilon)&\rightarrow A(1,1+\varepsilon)\\
    (r,e^{i\theta})&\mapsto(2+\varepsilon-r,-e^{-i\theta})
\end{align*}
$dF_l$ is the transition map we need. With respect to our trivializations on each side, this transition map becomes,
\begin{equation}\label{-2}
    (r,\theta,v)\mapsto\left(F_l(r,\theta),\left[\begin{matrix}
\cos (-2\theta) &-\sin(-2\theta)  \\
\sin (-2\theta) & \cos (-2\theta) \\
\end{matrix}\right]\cdot v\right).
\end{equation}
Based on this information, we construct a trivialization $\Phi_\Sigma:T\Sigma\rightarrow \Sigma\times\R^2$. Let $\Theta:\Sigma\rightarrow \text{U}(1)$ be a map whose winding number is -2 around each puncture that attaches $\Sigma_l$ and $-2g$ around the outer boundary. This outer boundary will be the attaching site for $T_{k+1}$. Also, let $\Theta$ be trivial away from those boundaries. Let a new trivialization $\Phi_D':TD\rightarrow D\times\R^2$ such that satisfies the following:
\begin{align*}
    \Phi_D'\circ\Phi_D^{-1}:(p,v)&\mapsto (p,\Theta(x)\cdot v).
\end{align*}
Then from \ref{-2}, $dF_l$ becomes the identity map via trivializations $\Phi_D'$ and $\Phi_{\Sigma_l}$. In other words, two trivializations are coherent, and thus one obtains a trivialization of $T\Sigma$. As we mentioned before, extend this trivialization of $T\Sigma$ to a trivialization of $\xi|_{\mathcal{E}}$ and name it $\Phi_{\mathcal{E}}$.

We now focus on attaching sites. Recall that we also used another coordinate $(s_j,u_j,t)$ on each gluing site $(1-\varepsilon,1)\times S^1\times S^1$ for $T_j$s. The Besse form was given by $\rho_j(s_j)du_j+dt$, and thus the symplectic form for $\xi$ restricted to the site is,
$$ \rho '_j(s_j)ds_j\wedge du_j.
$$
Note that according to the coordinate change on the attaching site of $T_j$, we have the following transition relation for $j=1,\dots,k$:
\begin{align*}\label{rel}
    \frac{\partial}{\partial x}&=-\cos(u_j)\cdot\frac{\partial}{\partial s_j}-\frac{\sin(u_j)}{s_j}\cdot\frac{\partial}{\partial u_j}\\
    \frac{\partial}{\partial y}&=\sin(u_j)\cdot\frac{\partial}{\partial s_j}-\frac{\cos(u_j)}{s_j}\cdot\frac{\partial}{\partial u_j}.
\end{align*}
For $j=k+1$ case, it turns out to be
\begin{align*}
    \frac{\partial}{\partial x}&=\cos(u_j)\cdot\frac{\partial}{\partial s_j}-\frac{\sin(u_j)}{s_j}\cdot\frac{\partial}{\partial u_j}\\
    \frac{\partial}{\partial y}&=\sin(u_j)\cdot\frac{\partial}{\partial s_j}+\frac{\cos(u_j)}{s_j}\cdot\frac{\partial}{\partial u_j}.
\end{align*}
Also, based on this relation, the following also holds:
\begin{equation}
    \Lambda(p)=\frac{\rho'(s_j)}{s_j}.
\end{equation}
Thus, according to the coordinate system $(s_j,u_j,t)$, the symplectic basis of $\xi$ near the $j$-th attaching site ($j=1,\dots,k$) consists of following elements:
\begin{align*}
   e_{j1}:=&-\frac{s_j\cos u_j}{\rho'(s_j)}\cdot\frac{\partial}{\partial s_j}-\frac{\sin u_j}{\rho'(s_j)}\cdot\frac{\partial}{\partial u_j}+\frac{\rho(s_j)\sin u_j }{\rho'(s_j)}\cdot\frac{\partial}{\partial t}\\
   e_{j2}:=&\sin u_j \cdot\frac{\partial}{\partial s_j}-\frac{\cos u_j}{s_j}\cdot\frac{\partial}{\partial u_j}+\frac{\rho(s_j)\cos u_j}{s_j}\cdot\frac{\partial}{\partial t}.
\end{align*}
For $j=k+1$ case, 
\begin{align*}
   e_{j1}:=&\frac{s_j\cos u_j}{\rho'(s_j)}\cdot\frac{\partial}{\partial s_j}-\frac{\sin u_j}{\rho'(s_j)}\cdot\frac{\partial}{\partial u_j}+\frac{\rho(s_j)\sin u_j }{\rho'(s_j)}\cdot\frac{\partial}{\partial t}\\
   e_{j2}:=&\sin u_j \cdot\frac{\partial}{\partial s_j}+\frac{\cos u_j}{s_j}\cdot\frac{\partial}{\partial u_j}-\frac{\rho(s_j)\cos u_j}{s_j}\cdot\frac{\partial}{\partial t}.
\end{align*}
The trivialization $\Phi_{\mathcal{E}}$ restricted to the $j$-th gluing site for $j=1,\dots,k$ uses the same basis above, but for $j=k+1$, it uses
\begin{align*}
    &\cos(2gu_{k+1})e_{(k+1) 1}+\sin (2gu_{k+1})e_{(k+1)2}  \\
    -&\sin(2gu_{k+1})e_{(k+1)1}+\cos(2gu_{k+1})e_{(k+1)2}.    
\end{align*}
   
For $\xi|_{T_j}$s, we trivialize the Besse structure as,
\begin{align*}
    &\Phi_{T_j}:\xi|_{T_j}\rightarrow T_j\times\R^2\\
    f_{j1}:=&\cos\theta_j\frac{\partial}{\partial r_j}-\frac{\sin\theta_j}{r_j}\cdot\frac{\partial}{\partial \theta_j}+\frac{\alpha_j r_j\sin \theta_j}{1+p_jr_j^2}\cdot\frac{\partial}{\partial \varphi_j}\mapsto (x,1,0)\\
    f_{j2}:=&\frac{1}{2}(1+p_jr_j^2)\left(\sin\theta_j\cdot\frac{\partial}{\partial r_j}+\frac{\cos\theta_j}{r_j}\cdot\frac{\partial}{\partial \theta_j}\right)-\frac{1}{2}\alpha_jr_j\cos\theta_j\frac{\partial}{\partial \varphi_j}\mapsto (x,0,1).
\end{align*}

Once we have our trivializations, we can now compute transition maps on each attaching site to get the first Chern class via obstruction theory. Since the Besse structure is a subbundle of the tangent bundle of $M$, its transition map is (the restriction of) the differential of the attaching map (\ref{gluingmap}). From this, one can obtain the transition map. Since the attaching cites are tori, all the transition maps can be written in the following form:
\begin{align*}
    T^2\times\R^2 &\rightarrow T^2\times\R^2\\
    (p,v)&\mapsto (p,m(p)\cdot v)
\end{align*}
Here, $p\mapsto m(p)$ is a smooth map from a torus to Sp(2), and we will always focus mainly on this matrix map while examining any transition maps. For $j=1,\dots,k$, such a matrix map can be described as the following product of matrices.
\begin{align*}
    \left[\begin{matrix}
\rho_j'(s_j) &0  \\
0 & s_j \\
\end{matrix}\right]
\cdot
\left[\begin{matrix}
\cos u_j &\sin u_j  \\
-\sin u_j & \cos u_j \\
\end{matrix}\right]
\cdot
\left[\begin{matrix}
-\frac{f'(r_j)}{s_j}&0  \\
0 & \frac{\alpha_j}{r_j(1+p_jr_j^2)} \\
\end{matrix}\right]
\cdot
\left[\begin{matrix}
\cos \theta_j &  \sin \theta_j \\
-\sin\theta_j & \cos \theta_j \\
\end{matrix}\right]
\cdot
\left[\begin{matrix}
1 &0  \\
0 & \frac{1}{2}(1+p_jr_j^2) \\
\end{matrix}\right]
\end{align*}
Note that (\ref{condi3}) was used during the computation. Now $s_j$ equals $f(r_j)$ and $u_j$ equals $-\alpha_j\theta_j-p_j\varphi_j$ according to our attaching map (\ref{gluingmap}). Also $\rho'_j(s_j),s_j,\frac{f'(r_j)}{s_j},\frac{\alpha_j}{r_j},$ and $1+p_jr_j^2$ are nonvanishing due to their definitions and conditions, including (\ref{condi1}),(\ref{small}), and (\ref{condi2}). Hence this map is homotopic along symplectic matrices to the following map.
\begin{align*}
\exp{i((\alpha_j-1)\theta_j+p_j\varphi_j)}.
\end{align*}
Here, we wrote the rotation matrix in the form of the exponential map for the sake of simplicity. For the $j=k+1$ case, we get the following matrix map.
\begin{align*}
    -&\left[\begin{matrix}
\cos 2gu_{k+1} &\sin 2gu_{k+1}  \\
-\sin 2gu_{k+1} & \cos 2gu_{k+1} \\
\end{matrix}\right]
\left[\begin{matrix}
\rho_{k+1}'(s_{k+1}) &0  \\
0 & s_{k+1} \\
\end{matrix}\right]
\cdot
\left[\begin{matrix}
\cos u_{k+1} &-\sin u_{k+1}  \\
\sin u_{k+1} & \cos u_{k+1} \\
\end{matrix}\right]\\
&\cdot
\left[\begin{matrix}
-\frac{f'(r_{k+1})}{s_{k+1}}&0  \\
0 & \frac{\alpha_{k+1}}{r_{k+1}(1+p_{k+1}r_{k+1}^2)} \\
\end{matrix}\right]
\cdot
\left[\begin{matrix}
\cos \theta_{k+1} & \sin \theta_{k+1} \\
-\sin\theta_{k+1} & \cos \theta_{k+1} \\
\end{matrix}\right]
\cdot
\left[\begin{matrix}
1 &0  \\
0 & \frac{1}{2}(1+p_{k+1}r_{k+1}^2) \\
\end{matrix}\right]
\end{align*}
Similarly, this map is homotopic to
\begin{equation*}
-\exp{i(((2g-1)\alpha_{k+1}-1)\theta_{k+1}+(2g-1)p_{k+1}\varphi_{k+1})}
\end{equation*}

At this point, we apply the obstruction theory \cite{SR}. According to the obstruction theory, transition maps of each 2-cells contain information about the Euler class (or, in other words, in our case, the first Chern class) of that plane bundle over the cellular complex. Once we know what each transition map looks like, up to homotopy, one can treat it as the following map:
\begin{align*}
    S^1\times\R^2&\rightarrow S^1\times\R^2\\
    (w,v)&\mapsto(w,w^r\cdot v).
\end{align*}
Here, $r$ is an integer, and this represents the winding number of the transition map. 2-cells in our case were each slices $C_j$ of $T_j$s with $\varphi_j$s being fixed. Its winding number for $j=1,\dots,k$ is then $\alpha_j-1$, and $(2g-1)\alpha_{k+1}$ for $j=k+1$. Then the cohomology class,
$$
     \begin{bmatrix}
        \alpha_1-1\\\alpha_2-1\\ \vdots \\\alpha_k-1\\ (2g-1)\alpha_{k+1}-1\\0
    \end{bmatrix}
    \in \Z^{k+2}/\,\,\emph{im}A_M \subseteq  H^2(M)
    $$
 becomes the first Chern class of the Besse structure by obstruction theory for the first Chern class (or, Euler class). (Note that the last entry becomes 0 as above because the last entry is for the 2-cell $\mathfrak{B}$, which was already included in $\mathcal{E}$.)
\begin{prop}\label{prev}
   The first Chern class of the Besse structure $\xi:=\ker\alpha$ of the Besse manifold $(M,\alpha)$ described by the Seifert invariant data $(0,g; (\alpha_1, \beta_1),\, \cdots \,(\alpha_{k+1},\beta_{k+1}))$ can be represented by,
   $$
     \begin{bmatrix}
        \alpha_1-1\\\alpha_2-1\\ \vdots \\\alpha_k-1\\ (2g-1)\alpha_{k+1}-1\\0
    \end{bmatrix}
    \in \Z^{k+2}/\,\,\emph{im}A_M \subseteq  H^2(M)
    $$
    where
    $$
    A_M=\begin{bmatrix}
    \alpha_1 & 0 & \cdots & 0& \beta_1 \\
0 & \alpha_2 & \cdots & 0 & \beta_2 \\
\vdots & \vdots & \ddots & \vdots & \vdots \\
0 & 0 & \cdots & \alpha_{k+1} & \beta_{k+1} \\
-1 & -1 & \cdots & -1 & 0 \\
\end{bmatrix}.
    $$
    (Note that the first component is for $\mathfrak{B}$, the 2-cell inside $\mathcal{E}$, so it becomes 0.)
    
   Especially, the Besse structure has trivial first Chern class $c_1(\xi)$ if and only if
    $$
     \begin{bmatrix}
        \alpha_1-1\\\alpha_2-1\\ \vdots \\\alpha_k-1\\ (2g-1)\alpha_{k+1}-1\\0
    \end{bmatrix}
    \in \emph{im} A_M.
    $$
\end{prop}

Now we rephrase the criterion for the first Chern class being trivial. This means we have integers $t_1,\dots,t_k,s$ such that,
$$
\begin{bmatrix}
        \alpha_1-1\\\alpha_2-1\\ \vdots \\\alpha_k-1\\ (2g-1)\alpha_{k+1}-1\\0
    \end{bmatrix}
    =
    \begin{bmatrix}
    \alpha_1 & 0 & \cdots & 0& 0& \beta_1 \\
0 & \alpha_2 & \cdots & 0 & 0 & \beta_2 \\
\vdots & \vdots & \ddots & \vdots &\vdots & \vdots \\
0 & 0 & \cdots & \alpha_k & 0 & \beta_k\\
0 & 0 & \cdots & 0 & \alpha_{k+1} & \beta_{k+1} \\
-1 & -1 & \cdots & -1&  -1 & 0 \\
\end{bmatrix}
\begin{bmatrix}
        t_1\\t_2\\ \vdots \\t_k\\t_{k+1}\\ s
    \end{bmatrix}
$$
The rational solution for this equation is
\begin{align}
    \begin{split}
     t_j&=\left(1-\frac{1}{\alpha_j}\right)+\frac{\beta_j}{\alpha_j}\cdot\frac{\chi}{e}\\
    s&=-\frac{\chi}{e}
\end{split}
\end{align}
where $\chi$ and $e$ are the orbifold Euler characteristic and Euler number, respectively, which we introduced in the preliminaries. 
\begin{prop}\label{crit}
    On the setting of Proposition \ref{prev}, the first Chern class of the Besse structure is trivial if and only if all of the following hold:
    \begin{enumerate}
        \item $\chi/e$ is an integer.
        \item For any $j$s, $\alpha_j|(\beta_j\cdot\chi/e-1)$.
    \end{enumerate}
\end{prop}

\begin{remark}\label{classification}
    Using the above criterion, we can classify Besse manifolds with vanishing first Chern class and non-negative orbifold Euler characteristic $\chi$. For $\chi>0$ case, lens spaces $L(p,p-1)$ for any positive integer $p$ and the contact structure induced from its smooth universal cover, the standard contact 3-manifold, are the only such cases. (Especially, this lens space becomes the standard contact 3-sphere when $p=1$.) In cases of vanishing $\chi$, Boothby-Wang bundles over a torus with positive Euler number are the result.  

    For the cases when $\chi$ is negative, it is difficult to obtain an explicit classification. For example, for any given $g>0$ and any $\alpha_j$s, we always have $c_1=0$ Besse manifolds with Seifert invariants $(2g-2,g;(\alpha_1,\alpha_1-1),\cdots,(\alpha_k,\alpha_k-1))$. (For $g=0$, this example still works with an additional condition stating that $\sum_j\frac{1}{\alpha_j}<k-2$.)
\end{remark}

\section{Invariance of symplectic homology}

We now prove our main theorem.

\begin{theorem}
    For a Besse manifold $(M,\xi=\text{ker }\alpha)$ with a given Seifert invariants
    
    \noindent $(b,g;(\alpha_1,\beta_1),(\alpha_2,\beta_2),\cdots,(\alpha_k,\beta_k))$,
    \begin{enumerate}[i]
        \item When the orbifold Euler characteristic is positive, the positive symplectic homology defined on its Stein filling becomes an invariant for such Besse manifolds.
        \item When the orbifold Euler characteristic is non-positive, the positive symplectic homology defined on the trivial Liouville cobordism $[0,1]\times M$ becomes an invariant for such Besse manifolds.
    \end{enumerate}
\end{theorem}
\begin{proof}
    For the first case, note that we have already classified this Besse contact manifold with $c_1=0$ (see \ref{classification}): Such Besse manifolds are $S^3$ with its standard tight contact structure and links of simple singularities. Both cases have a unique Stein filling with a trivial first Chern class. (\cite{G}, \cite{sim}) 

    For the non-positive case, note that a principal fiber in this case has infinite order in the fundamental group of $M$ \cite{Sei}. Since a suitable multiple of an exceptional fiber is homotopic to a principal fiber, one can see that all the Reeb orbits in this Besse manifold are noncontractible and hence the Besse contact structure is hypertight. Thus, we can define a positive symplectic homology over a trivial Liouville cobordism $[0,1]\times M$. (See \cite{cobo} for details.)
\end{proof}

\section{Robbin-Salamon indices of Reeb orbits}

In this section, we compute Robbin-Salamon indices of Reeb orbits in a Besse manifold $(M,\xi)$ with vanishing first Chern class. Again, suppose the Besse structure is based on Seifert invariants
$$
(b,g,(\alpha_1,\beta_1),\dots,(\alpha_k,\beta_k))
$$
or equivalently as before,
$$
(0,g,(\alpha_1,\beta_1),\dots,(\alpha_k,\beta_k),(\alpha_{k+1}=1,\beta_{k+1}=b)). 
$$
Since this Besse structure has vanishing first Chern class, it satisfies the criterion in Proposition \ref{crit}. Assume $W$ is a Stein filling of the Besse manifold $(M,\xi)$, and suppose that we have already attached the infinity cone part. (Later, we will find that any choice of Stein filling will not affect our result; see Remark \ref{indep}.)

Before the actual computation, we need a global trivialization for the given Besse structure to compute our indices under a coherent circumstance. Recall that each transition map on each attaching site for $T_j$s $(j=1,\dots,k)$ is homotopic to 
\begin{align*}
\exp{i((\alpha_j-1)\theta_j+p_j\varphi_j)}
\end{align*}
and for $j=k+1$ case is homotopic to 
\begin{equation*}
-\exp{i(((2g-1)\alpha_{k+1}-1)\theta_{k+1}+(2g-1)p_{k+1}\varphi_{k+1})}.
\end{equation*}
We now want to change the trivializations $\Phi_{\mathcal{E}}$ on the trivial bundle part $\mathcal{E}$ and $\Phi_{T_j}$ on each $T_j$ to make those transition maps into identity (up to homotopy). We only need to ``multiply" an appropriate map from $\mathcal{E}$ to Sp(2) for this, just as before. Note that by deformation-retracting $\mathcal{E}$ and Sp(2), each such map up to homotopy is 1-1 correspondence to each map up to homotopy from the product of a circle and the wedge of $k$ circles to a circle. Here, the $k$ circles for the wedge come from $k$ of those $k+1$ boundary components of $\Sigma$ and the other single circle in the product comes from a $S^1$-fiber of $\mathcal{E}$. Thus, up to homotopy, maps from $\mathcal{E}$ to Sp(2) are classified by the following invariants:
\begin{itemize}
    \item Winding numbers of the map along each of $k$ boundary components (with respect to orientations on each of the components). Without losing generality, suppose that we pick components where we attach $T_j$ for each $j=1,\dots,k$ and denote each winding number as $a_j$. 
    \item The winding number of the map along a single fiber with respect to the orientation coherent with that of the whole space. It will be denoted as $b$.
\end{itemize}
Let $\mathfrak{M}$ be one of such maps and denote its invariants $a_j$ $(j=1,\dots,k)$ and $b$ as above. Let our new trivialization $\Phi_{\mathcal{E}}':\xi|_{\mathcal{E}}\rightarrow\mathcal{E}\times\R^2$ on $\mathcal{E}$ such that the following holds:
$$
\Phi_{\mathcal{E}}'\circ\Phi_{\mathcal{E}}^{-1}:(p,v)\mapsto(p,\mathfrak{M}\cdot v).
$$
On the other hand, for some integer $c_j$ we also let $\Phi'_{T_j}:\xi|_{T_j}\rightarrow T_j\times\R^2$ be a new trivialization on each $T_j$ such that satisfies
$$
\Phi'_{T_j}\circ\Phi_{T_j}^{-1}:(p,v)\mapsto(p,\exp{ic_j\varphi}\cdot v).
$$ Regarding these new trivializations, one can compute the transition map up to homotopy for each $T_j$s similarly as before:
$$
\exp{i((\alpha_j-1)\theta_j+p_j\varphi_j+a_ju_j+bt-c_j\varphi_j)}.
$$
for $j=1,\dots,k$ and
$$
-\exp{i(((2g-1)\alpha_{k+1}-1)\theta_{k+1}+(2g-1)p_{k+1}\varphi_{k+1}+a_{k+1}u_{k+1}+bt-c_{k+1}\varphi_{k+1}))}
$$
for $j=k+1$ (where $a_{k+1}:=-(a_1+a_2+\dots+a_k)$, and $u_j=-\alpha_j\theta_j-p_j\varphi_j$ and $t=\beta_j\theta_j+q_j\varphi_j$ due to the attaching map).  The sufficient and necessary condition to make these maps into the identity map is,
\begin{align}
\begin{split}
    a_j&=\left(1-\frac{1}{\alpha_j}\right)+\frac{\beta_j}{\alpha_j}\cdot\frac{\chi}{e}\\
    b&=\frac{\chi}{e}\\
    c_j&=\frac{p_j}{\alpha_j}+\frac{1}{\alpha_j}\cdot\frac{\chi}{e}.    
\end{split}
\end{align}
Note that from the assumption that the first Chern class of $\xi$ is trivial, and our criterion Proposition \ref{crit}, all the above numbers become integers.

Once we have a trivialization, we can now compute the Robbin-Salamon indices of each Reeb orbit. Since the Besse structure exactly describes the Reeb flows, it is easy to see how the linearized Reeb flow goes with respect to our trivialization. We first look at a principal Reeb orbit. Note that any principal Reeb orbit is in $\mathcal{E}$, so one can compute its index only using coordinates on $\mathcal{E}$. Any such periodic orbit with multiplicity $m$ (i.e. an orbit that wraps around a principal fiber $m$ times) can be described as,
$$
t\in[0,1]\mapsto (x_0,y_0,\exp{2\pi mt})
$$
for some $(x_0,y_0)\in\Sigma$. Let $x_m$ be this orbit. Then, according to our trivialization, the linearized Reeb flow of time $t\in[0,1]$ along $x_m$ up to homotopy becomes,
$$
\begin{bmatrix}
    \cos(2\pi bmt) & -\sin(2\pi bmt) \\
\sin(2\pi bmt) & \cos(2\pi bmt) \\
\end{bmatrix}
$$
i.e., $bm=m\chi/e$ full rotations. Recall that for the path of rotation matrices $\mathcal{R}^{T}$ of $2\pi T$-rotation, the Robbin-Salamon index is

\begin{equation} \label{rot}
    \mu_{RS}(\mathcal{R}^{T})=
    \begin{cases*}
            2\left\lfloor T\right\rfloor+1& $T\in\R_{\geq0}\setminus\Z$ \\
            -2\left\lfloor-T\right\rfloor-1& $T\in\R_{\leq0}\setminus\Z$ \\
            2 T & $T\in\Z$\\
    \end{cases*}.
\end{equation}
Thus, the index in our case becomes,
$$
\mu(x_m)=2m\,\frac{\chi}{e}.
$$

For an exceptional Reeb orbit $x_{j,m}$ corresponds to a Seifert invariant $(\alpha_j,\beta_j)$ and multiplicity $m$, (\ref{action}) already describes the Reeb flow along the exceptional orbit. From this, one can compute the linearized Reeb flow of time $t\in[0,1]$ according to our new trivialization up to homotopy as, $$
\begin{bmatrix}
    \cos\left(2\pi \left(c_j-\frac{p_j}{\alpha_j}\right)mt\right) & -\sin\left(2\pi \left(c_j-\frac{p_j}{\alpha_j}\right)mt\right) \\
\sin\left(2\pi \left(c_j-\frac{p_j}{\alpha_j}\right)mt\right) & \cos\left(2\pi \left(c_j-\frac{p_j}{\alpha_j}\right)mt\right) \\
\end{bmatrix}
$$
i.e. $2\pi m(c_j-p_j/\alpha_j)=2\pi m(\chi/e)/\alpha_j$-rotation.
Note that the period of the exceptional orbit is $1/\alpha_j$ of that of the principal one, so this result makes sense. Hence, from (\ref{rot}), the Robbin-Salamon index for this exceptional orbit when $(m/\alpha_j)\cdot(\chi/e)$ is not integer is 
$$
\mu(x_{j,m})=
\begin{cases*}
    2\left\lfloor\frac{m}{\alpha_j}\cdot\frac{\chi}{e}\right\rfloor+1 & ($\chi>0$)\\
    -2\left\lfloor-\frac{m}{\alpha_j}\cdot\frac{\chi}{e}\right\rfloor-1 & ($\chi<0$)
\end{cases*}
$$
and if $(m/\alpha_j)\cdot(\chi/e)$ is an integer,
$$
\mu(x_{j,m})=2\frac{m}{\alpha_j}\cdot\frac{\chi}{e}.
$$

\begin{remark}\label{indep}
    One may have noticed that our choice of trivialization of $\xi$ was not unique nor canonical. As we explained at the beginning of this section, the maps from $\mathcal{E}$ to Sp(2) are up to homotopy represented by cohomology classes in $H^1(\mathcal{E};\Z)$. This cohomology group, as we have seen before, is generated by $k+2g$ generators, and we only used $k$ of them while constructing our trivialization. Changing the trivialization by adding or subtracting the leftover $2g$ generators is also possible and does not affect our computation. This also means that any choice of Stein fillings with vanishing first Chern class also does not affect the indices, unless it affects the homotopy class of each periodic orbit.
\end{remark}

\begin{remark}
    This index computation generalizes the result in \cite{Bries}, where the computation is done for Brieskorn manifolds, which are special cases of links of weighted homogeneous singularity (Example \ref{WH}). In fact, some of the $c_1(\xi)=0$ cases among Besse manifolds, such as links of $A_l,E_6,E_8$ singularities or circle bundles over a torus of Euler numbers $-1,-2,$ and $-3$ are actually Brieskorn manifolds.
\end{remark}

\section{Symplectic Homology of $\chi<0$ Cases}

Recall that for symplectic homologies, there is \emph{Morse-Bott spectral sequence} \cite{Bries}. Let $\{T_1,\dots,T_{r}\}$ be the set of minimal periods of all possible simple periodic Reeb orbits. Suppose that there is a minimal positive integer $n_1$ such that there exist positive integers $n_2,\dots,n_r$, and a real number $T_0$ satisfying $T_j=n_jT_0$ for all $j=1,\dots,r$. Let $C(p)$ be the set of Morse-Bott submanifolds with a period of $pT_0$, and let $s(B):=\mu_{RS}(B)-\frac{1}{2}\mathrm{dim}\,B/S^1$ for each Morse-Bott submanifold $B$, where $\mu_{RS}(B)$ is the Robbin-Salamon index of a periodic Reeb orbit in $B$.
\begin{equation}
    E_{pq}^1(SH^+)=
    \begin{cases*}
        \bigoplus_{B\in C(p)}H_{p+q-s(B)}(B;\mathbb{Z}/2\mathbb{Z}) & $p>0$ \\
        0 & $p\leq0$
    \end{cases*}
\end{equation}
This spectral sequence converges to positive symplectic homology.

Now, suppose we start from a Besse manifold, and suppose that this Besse manifold has a Stein filling with trivial first Chern class. For cases when $\chi$ is negative, one can observe that the degrees of orbits decrease as the action of it increases. If this degree  decreases fast, one would get an $E^1$-page with all the differentials being trivial, furthermore. And this can happen when $|\chi/e|$ becomes large enough compared to $\alpha_j$s. Hence, in such a case, we automatically get the symplectic homology, which is also independent of the choice of a Stein filling. (See Remark \ref{indep}.)

\begin{example}
    Let $M$ be a Besse manifold with Seifert invariants $(0,2;(2,1))$, and let $W$ be its simply-connected Stein filling with vanishing first Chern class. (Here, we just assume that such Stein filling exists.) Then $\chi/e=-5$, and thus this Besse structure is trivial as a plane bundle. From our index computation, the first page of the spectral sequence becomes the following.
    \begin{sseqdata}[ name = basic, xscale = 0.6, homological Serre grading ]
        \class(1,-6)
        \class(1,-5)
        \class(2,-13)
        \class(2,-12) \class(2,-12) \class(2,-12) \class(2,-12)
        \class(2,-11) \class(2,-11) \class(2,-11) \class(2,-11) 
        \class(2,-10)
        \class(3,-18) \class(3,-17)
    \end{sseqdata}
    \begin{center}
        \printpage[ name = basic, page = 1 ]
    \end{center}

    Since there cannot be non-trivial differentials on this page, the $E^2$-page remains the same, and so do the other pages. Thus, the positive symplectic homology of $M$ is then,
    \begin{equation*}
        SH^+_*(W;\mathbb{Z}/2\mathbb{Z})=
        \begin{cases*}
            (\mathbb{Z}/2\mathbb{Z})^{4\oplus} & $(*\leq -11, *\equiv 0,-1$ (mod 10)) \\
            \mathbb{Z}/2\mathbb{Z} & $(*\leq -4, *\equiv 2,5,6,9$ (mod 10))
        \end{cases*}
    \end{equation*}
    and this will be isomorphic to that of other Stein fillings.
\end{example}

\bibliographystyle{alpha}
\bibliography{ref}

@article{Wadsley,
    author = {A. Wadsley},
    title = {Geodesic Foliations by Circles},
    journal = {Journal of Differential Geometry},
    volume = {10},
    pages = {541--549},
    year = {1947}
}

@article{orbit,
    author = {D. Sullivan},
    title = {A counterexample to the periodic orbit conjecture},
    journal = {Publications math\'{e}matiques de l'I.H.\'{E}.S.},
    volume = {46},
    pages = {5--14},
    year = {1976}
}

@article{Abbondandolo_2022,
   title={Higher systolic inequalities for 3-dimensional contact manifolds},
   volume={9},
   journal={Journal de l'\'{E}cole polytechnique - Math\'{e}matiques},
   author={A. Abbondandolo and C. Lange and M. Mazzucchelli},
   year={2022},
    pages={807--851} }

@book{SR,
  author = {N. Steenrod},
  year = {1999},
  title = {The Topology of Fibre Bundles},
  publisher = {Princeton University Press}
}

@book{Geiges,
  author = {H. Geiges},
  year = {2008},
  title = {An introduction to contact topology},
  publisher = {Cambridge University Press},
  series = {Cambridge Studies in Advanced Mathematics, vol.109}
}

@article{ThurstonWinkelnkemper,
    author = {W. Thurston and H. Winkelnkemper},
    title = {On the Existence of Contact Forms},
    journal = {Proceedings of the American Mathematical Society},
    volume = {52},
    number = {1},
    pages = {345-347},
    year = {1975}
}

@article{Bries,
    author = {M. Kwon and O. van Koert},
    title = {Brieskorn manifolds in contact topology},
    journal = {Bulletin of the London Mathematical Society},
    volume = {48},
    number = {2},
    pages = {173-241},
    year = {2016}
}

@inbook{E, place={Cambridge}, series={London Mathematical Society Lecture Note Series}, title={Filling by holomorphic discs and its applications}, volume={2}, booktitle={Geometry of Low-Dimensional Manifolds: Symplectic Manifolds and Jones-Witten Theory}, publisher={Cambridge University Press}, author={Y. Eliashberg}, editor={S. Donaldson and C. Thomas}, year={1991}, pages={45--68}, collection={London Mathematical Society Lecture Note Series}}

@incollection{Abou,
    author = {M. Abouzaid},
    title = {Symplectic Cohomology and Viterbo's Theorem},
    booktitle = {Free loop spaces in geometry and topology},
    publisher = {European Mathematics Society},
    year = {2015}
}

@article{bourgeois2009symplectic,
  title={Symplectic homology, autonomous Hamiltonians, and Morse-Bott moduli spaces},
  author={Bourgeois, Fr{\'e}d{\'e}ric and Oancea, Alexandru},
  year={2009}, journal={Duke Mathematical Journal}, volume={146}, number={1}
}

@article{Ritter,
   title={Topological quantum field theory structure on symplectic cohomology},
   volume={6},
   number={2},
   journal={Journal of Topology},
   author={A. Ritter},
   year={2013},
   pages={391--489} }

@book{seifertt,
  author = {M. Jankins and W. N\'{e}umann},
  year = {1983},
  title = {Lectures on Seifert Manifolds},
  publisher = {Brandeis University, Waltham, MA},
  series = {Brandeis Lecture Notes 2}
}

@article{scott,
  title={The Geometries of 3-Manifolds},
  author={P. Scott},
  journal={Bull. London Math. Soc.},
  volume={15},
  number={5},
  pages={401--487},
  year={1983}
}

@article{SH1,
  title={Symplectic homology. \text{II}. A general construction},
  author={K. Cieliebak and A. Floer and H. Hofer},
  journal={Math.Z.},
  volume={218},
  pages={103--122},
  year={1995}
}

@article{SH0,
    author = {A. Floer and H. Hofer},
    title = {Symplectic homology. \text{I}. Open sets in {$\C^n$}},
    journal = {Mathematische Zeitschrift},
    year = {1994},
    volume = {215},
    pages ={37--88}
}

@article{SH2,
  title={Functors and computations in Floer homology with applications. \text{I}.},
  author={C. Viterbo},
  journal={Geom. Funct. Anal.},
  volume={9},
  pages={985--1033},
  year={1999}
}

@book{AD,
  author = {M. Audin and M. Damian},
  year = {2014},
  title = {Morse Theory and Floer Homology},
  publisher = {Springer}  
}

@article{cobo,
   title={Symplectic homology and the Eilenberg-Steenrod
axioms},
   volume={18},
   number={4},
   journal={Algebraic and Geometric Topology},
   author={K. Cieliebak and A. Oancea},
   year={2018}, pages={1953--2130} }

@book{sasaki,
  author = {C. P. Boyer and K. Galicki},
  year = {2008},
  title = {Sasakian geometry},
  publisher = {Oxford University Press, Oxford},
  series = {Oxford Mathematical Monographs}
  
}

@phdthesis{mcleanstein,
    author = {M. McLean},
    title = {The Symplectic Topology of Stein Manifolds},
    school = {University of Cambridge},
    year = {2008}
}

@article{Fauck,
    author = {A. Fauck},
    title = {Rabinowitz Floer Homology on Brieskorn Spheres},
    journal = {International Mathematics Research Notices},
    volume = {2015},
    issue = {14},
    Pages = {5874--5906},
    year = {2015}
}

@phdthesis{Guttphd,
    author = {J. Gutt},
    title = {On the minimal number of periodic Reeb orbits on a contact manifold},
    school = {Universit\'{e} de Strasbourg},
    year = {2014}
}

@article{Guttpesh,
    author = {J. Gutt},
    title = {The positive equivariant symplectic homology as an invariant for some contact manifolds},
    journal = {Journal of symplectic geometry},
    year = {2017},
    volume = {15},
    number = {4},
    pages = {1019--1069}
}

@article{U1,
    author = {I. Ustilovsky},
    title = {Infinitely many contact structures on {$S^{4m+1}$}},
    journal = {Internat. Math. Res. Notices},
    number = {14},
    pages = {781--791},
    year = {1999}
}

@article{Zhou_2021,
   title={$(\mathbb{RP}^{2n-1},\xi_{\text{std}})$ is not exactly fillable for
$n\neq 2^k$},
   volume={25},
   number={6},
   journal={Geometry and Topology},
   author={Z. Zhou},
   year={2021},
   pages={3013--3052} }

@article{Zhou_2024,
   title={On fillings of contact links of quotient singularities},
   volume={17},
   number={1},
   journal={Journal of Topology},
   author={Z. Zhou},
   year={2024},}

@article{Benedetti_2020,
   title={Invariance of symplectic cohomology and twisted cotangent bundles over surfaces},
   volume={31},
   number={09},
   journal={International Journal of Mathematics},
   author={G. Benedetti and A. Ritter},
   year={2020},
   pages={2050070} }

@InProceedings{biased,
   title={A biased view of symplectic cohomology},
   booktitle={Current Developments in Mathematics, 2006},
   author={P. Seidel},
   year={2008},
    publisher={Int. Press},
   pages={211--253} }

@article{sim,
author = {H. Ohta and K. Ono},
title = {{Simple Singularities and Symplectic Fillings}},
volume = {69},
journal = {Journal of Differential Geometry},
number = {1},
publisher = {Lehigh University},
pages = {001 -- 42},
year = {2005},
}

@book{CE,
  author    = {K. Cieliebak and Y. Eliashberg},
  title     = {From Stein to Weinstein and back: symplectic geometry of affine complex manifolds},
  publisher = {American Mathematical Society},
  series = {Colloquium Publications},
  volume = {59},
  year      = {2012},
}

@book{orlik,
  author = {P. Orlik},
  year = {1972},
  title = {Seifert Manifolds},
  publisher = {Springer Berlin, Heidelberg},
  series = {Lecture Notes in Mathematics}
}

@article{BW,
  title={On contact manifolds},
  author={W. Boothby and H. Wang},
  journal={Ann. of Math.},
  volume={68},
  number={2},
  pages={721--734},
  year={1958}
}

@article{Kegel,
  title={A Boothby-Wang theorem for Besse contact manifolds},
  author={M. Kegel and C. Lange},
  journal={Arnold Math. J.},
  volume={7},
  pages={225--241},
  year={2021},
}

@article{G,
  title={Pseudoholomorphic curves in symplectic manifolds},
  author={M. Gromov},
  journal={Invent. Math.},
  volume={4},
  number={2},
  pages={307--347},
  year={1985},
  publisher={International Press}
}

@article{CZ,
  title={Morse Theory for periodic solutions of Hamiltonian systems and Maslov index},
  author={D. Salamon and E. Zehnder},
  journal={Comm. Pure Appl. Math},
  volume={45},
  pages={1303--1360},
  year={1992}
}

@article{RS,
  title={The Maslov index for paths},
  author={J. Robbin and D. Salamon},
  journal={Topology},
  volume={32},
  pages={827--844},
  year={1993},
  number={4}
}

@article{simE,
title = {Symplectic fillings of the link of simple elliptic singularities},
author = {H. Ohta and K. Ono},
pages = {183--205},
volume = {2003},
number = {565},
journal = {Journal für die reine und angewandte Mathematik},
year = {2003},
}

@article{Sei,
    author = {H. Seifert},
    title = {Topology of 3-dimensional fibered spaces},
    journal = {Acta Mathematica},
    volume = {60},
    pages = {147--288},
    year = {1933}
}

@article{CR,
    author = {L. Lempert},
    title = {On three-dimensional Cauchy-Riemann manifolds},
    journal = {Journal of the american mathematical society},
    year = {1992},
    volume = {5},
    number = {4}
}

\end{document}